\documentclass[12pt]{article}
\usepackage{amsmath} %
\usepackage[dvips]{color}
\usepackage{subfigure}
\usepackage{graphicx}
\definecolor{dgreen}{rgb}{0,.8,.3}
\definecolor{lblue}{rgb}{.2,.3,.7}

\newcommand{\bull}{\vrule height 1.8ex width 1.0ex depth 0ex}

\begin{document}
\begin{center}

{\bf\Large Sparse Index Tracking Based On $L_{1/2}$ Model And Algorithm } \vspace{0.5in}

 Xu Fengmin \footnote{Xu Fengmin. School of Mathematics and Statistics,
 Xi'an Jiaotong University, Xi'an, 7 10049, P. R. China. E-mail: fengminxu@mail.xjtu.edu.cn. This work is supported
 by China NSFC projects No. 11101325 and Open Fund of Key Laboratory of Precision Navigation and Technology, National Time Service Center, CAS. }
\qquad Zongben Xu\footnote{Zongben Xu. School of Science, Xi'an Jiaotong University, Xi'an, 710049, China.~E-mail: zbxu@mail.xjtu.edu.cn.~Supported by the National 973 Program of China(No.2007CB311002),
the China NSFC projects(No.70531030).}
\qquad Honggang Xue\footnote{Honggng Xue. School of Economics and Finance,
 Xi'an Jiaotong University, Xi'an, 710049, China.~E-mail:
xhg@mail.xjtu.edu.cn.~Supported by the Ministry of Education
Humanities Social Science(No.09XJAZH005), the Ministry of Education
new century elitist supports plan(No.NCET-10-0646) }

\vspace{0.3in}
\end{center}

%
%
%
%

\noindent {\bf Abstract.}  Recently, $L_1$ regularization have been attracted extensive attention and successfully applied in mean-variance portfolio selection for promoting out-of-sample properties and decreasing transaction costs. However, $L_1$ regularization approach is ineffective in promoting sparsity and selecting regularization parameter on index tracking with the budget and no-short selling constraints, since the 1-norm of the asset weights will have a constant value of one. Our recent research on $L_{1/2}$ regularization has found that the half thresholding algorithm with optimal regularization parameter setting strategy is the fast solver of $L_{1/2}$ regularization, which can provide the more sparse solution. In this paper we apply $L_{1/2}$ regularization method to stock index tracking and establish a new sparse index tracking model. A hybrid half thresholding algorithm is proposed for solving the model. Empirical tests of model and algorithm are carried out on the eight data sets from OR-library. The optimal tracking portfolio obtained from the new model and algorithm has lower out-of-sample prediction error and consistency both in-sample and out-of-sample. Moreover, since the automatic regularization parameters are selected for the fixed number of optimal portfolio, our algorithm is a fast solver, especially for the large scale problem.

 \vspace{0.1in}

\vspace{0.1in} \noindent {\bf AMS Subject Classifications:}

\vspace{0.2in}
\section{Introduction}

Stock index derivatives, such as index funds, index futures, index options
etc, have developed very rapidly and become important tools in investment
and risk management of global financial markets, especially it shows the
better effects to stabilize the stock market in global fiance crisis. Index
tracking (e.g., index replication) plays a core role in product design and
risk management of index derivatives. It consists in construction of a
tracking portfolio whose behavior is as similar as possible to target index
during a predefined period.

Broadly speaking, two different strategies can be used to track a given
stock market index: the full replication and the non-full replication. The
full replication consists in purchasing all constituent stocks of a given
index. In practice, this strategy need high transaction costs. An
alternative way is the non-full replication, which include the stratified
sampling replication and the optimal replication. Since the selection of the
stocks in stratified sampling replication depends on the manager's
experience, so the tracking portfolio is non-optimal, thus we focus on the
optimal replication method in this paper. The optimal replication aims to
find the portfolio that minimizes the tracking error by investing in only a
subset of the assets using optimization method. This strategy involves much
lower transaction costs, and can achieve acceptable tracking errors in
principle.

Different quantitative methods have been proposed to tackle such an
optimization problem. Roll establishes optimal index tracking models and
proposed a mean-variance analysis of index tracking on Markowitz's earlier
study \cite{Roll}. Tabata and Takeda discuss the index fund management based
on mean-variance model \cite{Tabata}. Buckley and Korn apply optimal impulse
control techniques to the index tracking problem with fixed and proportional
transaction costs \cite{Korn}. Rudolf et al. propose several piecewise
linear measures of the tracking error, and solve the problem by means of
linear programming \cite{Rudolf}. Alexander proposes the construction of
tracking portfolios by analyzing the coincidental structure between the time
series of each of the assets and the time series of the tracked index \cite
{Alexander}. Ammann and Zimmermann investigate the relationship between
several statistical measures of tracking error and asset allocation
restrictions based on admissible weight ranges \cite{Ammann}. Gilli and
Kellezi propose the use of the threshold accepting heuristic to solve the
problem, including cardinality restrictions and transaction costs \cite%
{Gilli}. Beasley et al. address the index tracking problem using
evolutionary heuristics with real-valued chromosome representations \cite%
{Beasley2003} . Lobo et al. investigate the portfolio optimization problem
with transaction costs, which they address by means of a heuristic
relaxation method that consists in solving a small number of convex
optimization problems using fixed transaction costs \cite{lobo}. Torrubiano
and Beasley present a nonlinear mixed-integer optimal model and a
corresponding algorithm for index tracking \cite{Beasley}. Torrubian et al.
design a hybrid strategy that combines an evolutionary algorithm with
quadratic programming to yield the optimal tracking portfolio that invests
only in the selected assets \cite{Torr}.

On the other hand, statistical regularization methods have been successfully
applied in mean-variance portfolio selection in order to promote the
identification of sparse portfolios with good out-of-sample properties and
low transaction costs \cite{DeMiguel, Brodie, Fanjianqing}. DeMiguel et al
focuses on the effect of the constraints on the covariance regularization, a
technique extension of the result in Jagannathan and Ma \cite{Ma}. Brodie et
al emphasize on the sparsity of the portfolio allocation and the
optimization algorithms by using the LASSO ($L_{1}$ regularization \cite{43}%
), they also noted that the idea using $L_{1}$ regularization can be used to
solve the index tracking problem with short selling constraints. Prominent
contribution of Fan Jianqing et al is to provide mathematical insights to
the utility approximations with the gross-exposure constraint. These
proposed approaches rely on imposing upper bounds on the 2-norm of the
portfolio weights as suggested by the ridge regression($L_{2}$
regularization \cite{Hoerl}), or on the 1-norm using $L_{1}$ regularization
approach. Empirical results in a mean-variance framework support the use of
the $L_{1}$ regularization method when short selling is allowed. However,
the LASSO approach is ineffective in promoting sparsity in presence of the
budget and no-short selling constraints.

Consider the index tracking problem, the budget and no-short selling
constraints is essential. If we use the $L_{1}$ regularization to deal with
the index tracking problem, there will be some defects. First, the $L_{1}$
regularization can't provide the more sparse optimal solution since the $1$%
-norm of the asset weights will have a constant value of one; Second, the
selection of regularization parameter is a hard problem for $L_{1}$
regularization since the number of the optimal tracking portfolio is fixed;
Finally, the optimization strategy to deal with the constrained $L_{1}$
regularization is to use the penalty function method, the penalty factor is
more difficult to select.

Fortunately, our recent studies on $L_{1/2}$ regularization have found that $
L_{1/2}$ regularization can overcome these defects of $L_{1}$ regularization
\cite{52,53,54}. The reasons are as follows. Firstly, using $L_{1/2}$
regularization get the more sparse tracking portfolio than $L_{1}$
regularization \cite{Fernholz}, that is we can use the least stocks to track
the target index by controlling the turnover; Secondly, though $L_{1/2}$
regularization is a nonconvex, non-smooth and non-Lipschitz optimization
problem, we derive the fast and effective half thresholdig algorithm for
solution of $L_{1/2}$ regularization, especially for large-scale problems
\cite{54}. Finally, For decreasing transaction costs and easy to manage
portfolio, managers often request a sparse tracking portfolio with fixed $K$
stocks to track the object index when index has a large number of
constituents. For $K$-sparsity index tracking problem, the regularized
parameter of half thresholding algorithm can automatic correct to
appropriate value whatever the initial condition is.


Base on the above analysis, The main work in this paper is to design a
sparse index tracking model and algorithm by introducing $L_{1/2}$
regularization. Different to focus on finding the portfolio that is optimal
using as inputs the recent historical evolution of the assets, we are
interested in the future tracking performance of the portfolio. In section
2, we briefly review the index tracking model in \cite{Torr} and $L_{1/2}$
regularization with half thresholding algorithm. A sparse index tracking
model with a hybrid half thresholding algorithm based on $L_{1/2}$
regularization and half thresholding algorithm is derived in section 3.
Empirical comparisons are conducted in section 4. The data is partitioned
into training data and testing data, the training data is used to construct
the optimal tracking portfolio investing in a subset of the index assets.
The performance of this tracking portfolio is then evaluated not only on the
in-sample data, but also on the out-of-sample data. So the optimal tracking
portfolio of sparse index tracking model that are suboptimal on the training
data can have a better out-of-sample performance on the test data.
Meanwhile, we define the consistence indicator to discuss the performance of
the new model both in-sample and the out-of-sample, the empirical comparison
in section 4 will illustrate these results. We conclude the paper in section
5.


\section{Preliminaries}

\setcounter{equation}{0} To give a precise formulation of the sparse index
tracking model and algorithm, we first review the index tracking model from
the regression viewpoint \cite{Torr}, then we provide a general account of
the $L_{1/2}$ regularization and the half thresholding algorithm which serve
the basis of the new model and algorithm.


\subsection{Index tracing problem}

In this subsection we give the index tracking model by introducing the
tracking error which is treated as the objective function and constraints of
the index tracking problem from the the regression viewpoint.

The tracking error have many different definitions, consequently, different
tracking portfolio models are introduced, see \cite{Ammann, lobo, Shapcott,
Beasley2003, Korn, Rudolf}. Most of them introduce the definition of
tracking error based either on correlations between the returns of tracking
portfolio and the index or on estimates of the variance between the returns
of the index and the returns of the tracking portfolio \cite{Korn, Rudolf}.
However, Beasley et al. argue against the use of variance as a measure of
tracking error because the tracking error would be zero while the difference
between the return of the index and the tracking portfolio is constant \cite%
{Beasley2003}. This is the undesirable result because it does not take into
account the tracking bias. Beasley et al. give a new definition of tracking
error, that is, the square of mean squared error to measure the difference
between the return of the target index and the tracking portfolio \cite%
{Beasley2003}, this definition of the tracking error takes into account the
bias of the tracking portfolio. Consequently, we adopt this definition of
tracking error in this paper.

Let $P_{it}$ be the time series of stock prices for the $N$ stocks that are
included in the given stock market index whose evolution we wish to
replicate. Let $I(t)$ be the time series of this index. All time series are
defined for equally spaced intervals $t=1,2,...,T$. Under fixed mixture
strategy, the tracking error is defined
\begin{equation}
\begin{array}{cc}
TE=\frac{1}{T}\sum\limits_{t=1}^{T}(\sum%
\limits_{i=1}^{N}w_{i}r_{it}-R_{t}^{I})^{2}, &
\end{array}
\label{21}
\end{equation}%
where

$r_{it}:$ \ \ the return rate of stock $i$ at time $t$ during single period,
that is
\begin{equation}  \label{22}
\begin{array}{cc}
r_{it}=\frac{P_{it+1}-P_{it}}{P_{it}};\ \ i=1, \cdots, N, \ \ t=1,2,\cdots,T.
&
\end{array}%
\end{equation}

$R_{t}^{I}:$\ \ the return rate of the target index at time $t$ during
single period, that is
\begin{equation}  \label{23}
\begin{array}{cc}
R_{t}^{I}=\frac{I_{t+1}-I_{t}}{I_{t}};\ \ t=1,2,\cdots,T. &
\end{array}%
\end{equation}

$w_{i}:$\ \ the weights of stock $i$.\newline

Let
\begin{equation*}
R^{I}=(R_{1}^{I},R_{2}^{I}\cdots ,R_{T}^{I})^{T}\in R^{T\times 1}
\end{equation*}%
the column vector of the index return rate, and
\begin{equation*}
R=(R_{1},R_{2},\cdots ,R_{N})\in R^{T\times N}=\left(
\begin{array}{cccl}
r_{11}\, & r_{21}\, & \cdots \, & r_{1N} \\
r_{12}\, & r_{22}\, & \cdots \, & r_{1N} \\
\vdots \, & \cdots \, & \cdots \, & \vdots  \\
r_{1T}\, & r_{2T}\, & \cdots \, & r_{NT}%
\end{array}%
\right)
\end{equation*}%
where $R_{i}=(r_{i1},r_{i2},\cdots ,r_{iT})^{T}$ is the column vector of
return rate of the stock $i$ , $i=1,2,\cdots ,N$. $R$ is the matrix of all
stock's return rate. Let $w=(w_{1},w_{2},\cdots ,w_{N})^{T}$ be the $N\times
1$ column vector of the stock weights, the tracking error (\ref{21}) can be
replaced as
\begin{equation}
\begin{array}{cc}
TE=\frac{1}{T}\Vert Rw-R^{I}\Vert _{2}^{2}. &
\end{array}
\label{24}
\end{equation}

The aim of index tracking problem is to find optimal tracking portfolio by
minimizing the tracking error (\ref{24}) under some constraints.

The first constraints of index tracking model is the budget constraints
\begin{equation*}
\sum\limits_{i=1}^{N}w_{i}=1,
\end{equation*}%
it ensures that all the capital is invested in the tracking portfolio. The
second item is the lower and upper bound constraints
\begin{equation}
\begin{array}{cc}
\eta _{i}Z_{i}\leq w_{i}\leq Z_{i}\delta _{i},\ \ i=1,\cdots ,N. &  \\
&
\end{array}
\label{25}
\end{equation}%
The aim of setting lower bounds of the investment ratio $w_{i}$ is to avoid
small invest volume, and setting upper bounds is to control risk. The third
item is the cardinality constraints
\begin{equation}
\sum\limits_{i=1}^{N}Z_{i}=K,~~~~Z_{i}=0\ \ or\ 1,~~~~i=1,2,\cdots ,N,
\label{26}
\end{equation}%
where $K$ is the number of the stocks included in a tracking portfolio which
K is a given positive integer. (\ref{26}) reflects that if asset $i$ is not
included in the tracking portfolio then $Z_{i}=0$, otherwise, $w_{i}=0$ by (%
\ref{25}).

Based on the definition of tracking error and the constraints, the basic
index tracking model is described as \cite{Torr}
\begin{equation}  \label{27}
\begin{split}
\text{min} & ~~\frac{1}{T}\|Rw-R^I\|^2_2 \\
\text{s.t.}& ~~e^Tw=1 \\
& ~~\eta_{i}Z_{i}\leq w_{i}\leq Z_{i}\delta_{i} \\
& ~~\sum\limits_{i=1}^{N}Z_{i}=K \\
& ~~Z_{i}=0 \ \ or \ 1,~~i=1,2,\cdots,N.
\end{split}%
\end{equation}

The above index tracking problem is ERM model, the optimal tracking
portfolio obtained (\ref{27}) has the minimal in-sample tracking error, but
we can't know the performance of out-of-sample error. The model (\ref{27})
is hard to solve since the cardinally constraint ($\sum%
\limits_{i=1}^{N}Z_{i}=K$) is discrete and therefore highly nonlinear, many
different optimization technique have been proposed to tackle such a hybrid
nonlinear integer programming, see \cite{Beasley2003,Beasley,Torr}.

Different to the optimization technique, we hope to give the new sparse
index tracking model which is easily solved by introducing the $L_{1/2}$
regularization. The new model can generate the sparse solution with good
tracking performance both in-sample and out-of-sample. For constructing our
model and solving it efficiently, we review the $L_{1/2}$ regularization and
half thresholding algorithm in the next subsection.


\subsection{$L_{1/2}$ regularization}

In this subsection, we briefly introduce the $L_{1/2}$ regularization and
the half thresholding algorithm \cite{52,53,54}, and explain why we use $%
L_{1/2}$ regularization for solving the index tracking problem.

$L_{1/2}$ regularization is one of the statistical regularization methods
that is used to solve the sparse problem, which aims to find sparse solution
of a representation or an equation. Typically, the sparsity problems include
those of variable selection \cite{43}, visual coding \cite{39,30}, graphical
modeling \cite{32}, error correction \cite{8}, matrix completion \cite{7}
and compressed sensing \cite{53, 11, 9, 19}.

$L_{1/2}$ regularization can be modeled as the following optimization problem

\begin{equation}  \label{288}
\begin{array}{cc}
\min\limits_{x\in R^{N}} \|Ax-b\|^{2}+\lambda \|x\|^{1/2}_{1/2}, &
\end{array}%
\end{equation}
where $A\in R^{M\times N},\ \ x=(x_{1},\cdots, x_{N})^{T}\in R^{N}$, $%
\|x\|_{1/2}^{1/2}=\sum\limits_{i=1}^{n}|x_{i}|^{1/2}.$ $\lambda$ is the
regularization parameter which control the sparsity of optimal solution (\ref%
{288}).

In general, $L_{0}$ regularization and $L_{1}$ regularization are also the
efficient method for solving the sparsy problem. $L_{0}$ regularization is
\begin{equation}
\begin{array}{cc}
\min\limits_{x\in R^{N}}\Vert Ax-b\Vert ^{2}+\lambda \Vert x\Vert _{0} &
\end{array}
\label{28}
\end{equation}%
where $\Vert x\Vert _{0}$ means the number of nonzero components in $x$ . $%
L_{1}$ regularization is
\begin{equation}
\begin{array}{cc}
\min\limits_{x\in R^{N}}\Vert Ax-b\Vert ^{2}+\lambda \Vert x\Vert _{1} &
\end{array}
\label{29}
\end{equation}%
where $\Vert x\Vert _{1}$ means the 1-norm of $x.$

Unfortunately, $L_{0}$ regularization (\ref{28}) is NP-hard and hardly
tractable when $x$ is large. $L_{1}$ regularization (\ref{29}), known as the
Lasso, have been introduced in the ninthes by Tibshirani \cite{43} and it
has also been independently proposed by Chen et al \cite{399} as the basis
pursuit denoising problem. $L_{1}$ regularization is the convex optimization
problem and have the analytic solution. However, although the $L_{1}$
regularization provides the best convex approximation to the $L_{0}$
regularization and is computationally efficient, the $L_{1}$ regularization
cannot handle collinearity and may result in inconsistent selection and
introduce extra bias in estimation. One valid improved method is to use $%
L_{1/2}$ regularization, the $L_{1/2}$ regularization can generate the more
sparse solutions than $L_{1}$ regularization.

Though $L_{1/2}$ regularization (\ref{288}) leads to a nonconvex, non-smooth
and non-Lipschitz optimization problem, our recent studieds \cite{52, 53, 54}
dissolved this problem. Through justifying the existence of the resolvent of
gradient of penalty, and looking for its analytic expression, we derived an
iterative half thresholding algorithm \cite{54} for fast solution of $%
l_{1/2} $ regularization. We prove an alternative feature theorem on
solutions of $L_{1/2}$ regularization, based on which a thresholding
representation of $L_{1/2}$ regularization is given and a novel
regularization parameter setting strategy is suggested. We verify the
convergence of the iterative half thresholding algorithm and provide a
series of experiments and applications to assess performance of the
algorithm.

The half thresholding algorithm can be described as
\begin{equation}
x_{n+1}=H_{\lambda _{n}\mu _{n},1/2}(x_{n}+\mu _{n}A^{T}(b-Ax_{n})).
\label{66}
\end{equation}%
where $H_{\lambda \mu ,1/2}(x)=(h_{\lambda \mu ,1/2}(x_{1}),h_{\lambda \mu
,1/2}(x_{2}),\cdots ,h_{\lambda \mu ,1/2}(x_{N}))$ is the half thresholding
operator, for $i=1,\cdots ,N$
\begin{equation*}
h_{\lambda \mu ,1/2}(x_{i})=\left\{
\begin{array}{cc}
\frac{2}{3}|x_{i}|\left( 1+\mathrm{cos}\left( \frac{2\mathrm{\pi }}{3}-\frac{%
2\varphi _{\lambda }(x_{i})}{3}\right) \right) , & \,|x_{i}|>\frac{\sqrt[3]{%
54}}{4}(\lambda \mu )^{\frac{2}{3}} \\
0, & \,otherwise%
\end{array}%
\right.
\end{equation*}%
and
\begin{equation}
\cos \varphi _{\lambda }(x_{i})=\frac{\lambda }{8}\left( \frac{|x_{i}|}{3}%
\right) ^{-\frac{3}{2}}.
\end{equation}%
It is known that the quantity of solutions of a regularization problem
depends seriously on the setting of regularization parameter $\lambda .$ The
selection of proper regularization parameters is, however, a very hard
problem. However, If the solutions of problem (\ref{288}) are $K$-sparsity,
the parameters can be found by formulating an optimality condition on
regularization. Let $B_{\mu }(x_{n})=I+\mu A^{T}(b-Ax_{n}),$ and $[B_{\mu
}(x_{n})]_{k+1}$ is the $(K+1)$-th largest component of $B_{\mu }(x_{n})$ in
magnitude, the parameters can be adopted with
\begin{equation*}
\lambda _{n}=\frac{\sqrt{96}}{9\mu _{0}}\big |[B_{\mu _{0}}(x_{n})]_{k+1}%
\big|^{\frac{3}{2}},
\end{equation*}%
and the constant $\mu _{n}=\mu _{0}>0$.

In fact, the problem of tracking a financial index using only a subset of
stocks can be regarded as the sparsity problem. If the number of stocks
included in the tracking portfolio is fixed, the problem of selecting the
optimal $K$ stocks is $K$-sparsity problem. and then can be handled by the $%
L_{1/2}$ regularization. Different to use $L_{1}$ regularization, the index
tracking problem using $L_{1/2}$ regularization can provide more sparse
tracking portfolio. Furthermore, the $L_{1/2}$ regularization has the fast
and efficient algorithm with the better method of regularization parameter
selection for the $K$-sparsity problem. Hence, the new formulation of index
tracking problem using $L_{1/2}$ regularization are presented in next
section.


\section{{\protect\large A $L_{1/2}$ regularization based model}}

\setcounter{equation}{0}

In this section, we propose a new sparse index tracking model and the hybrid
half thresholding algorithm based on the analysis of section 2.

\subsection{The sparse index tracing model}

Considering the constraints of the index tracking model (\ref{27})
introduced in section 2, Let
\begin{equation*}
\Omega_{1}=\{w|Z_{i}\eta_{i}\leq w_{i}\leq Z_{i}\delta_{i},\
\sum\limits_{i=1}^{N}Z_{i}=K,\ Z_{i}=0 \ \ or \ 1,\ \ i=1,\cdots, N,\}
\end{equation*}
the constraints $\sum\limits_{i=1}^{N}Z_{i}=K,\ Z_{i}=0 \ \ or \ 1,\ \
i=1,\cdots, N$ means the number of nonzero components of the optimal
tracking weight $w$ is $K,$ and the tracking weight $w_{i}=0 $ if $Z_{i}=0$
or $\eta_i\leq w_{i}\leq \delta_i$ if $Z_i=1$ for $i=1,\cdots,N.$ We also
notice that $\|w\|_0$ means the number of nonzero components of the optimal
tracking weight $w,$ then let
\begin{equation*}
\Omega_{2}=\{w| \|w\|_0=K, w_{i}=0 ~~ or ~~\eta_i\leq w_{i}\leq \delta_i,
~~i=1,\cdots,N\}.
\end{equation*}
Let $Supp(w)$ be the support set of $w,$ i.e., $Supp(w)=\{i|w_i\neq 0\},$ a
sparse, stable index tracking model is obtained by adopting a regularization
procedure, that is, the constraint set $\Omega_{1}$ is replaced by $%
\Omega_{2}$, so the sparse index tracking model based on the $l_{0}$
regularization is
\begin{equation}  \label{31}
\begin{split}
\text{min} & ~~\frac{1}{T}\|Rw-R^I\|^2_2 \\
\text{s.t.}& ~~e^Tw=1 \\
& ~~\|w\|_0=K \\
& ~~\eta_i\leq w_{i}\leq \delta_i, \ \ i\in Supp(w) \\
&~~w_{i}=0, \ \ i\notin Supp(w)
\end{split}%
\end{equation}

Based on the analysis of the subsection 2.2, it is better that we use $
\|w\|_{1/2}^{1/2}=K$ to substitute the constraint $\|w\|_{0}=K$, and we omit
the coefficient $\frac{1}{T}$ of the objective function. Then a new sparse
index tracking model based on $L_{1/2}$ regularization is proposed as
\begin{equation}  \label{32}
\begin{split}
\text{min} & ~~\|Rw-R^I\|^2_2 \\
\text{s.t.}& ~~e^Tw=1 \\
& ~~\|w\|_{1/2}^{1/2}=K \\
& ~~\eta_i\leq w_{i}\leq \delta_i, \ \ i\in Supp(w) \\
&~~w_{i}=0, \ \ i\notin Supp(w)
\end{split}%
\end{equation}

In order to solve this model efficiently, we penalize the constraint $%
\|w\|_{1/2}^{1/2}=K$ to the objective function using penalty function
method, then an equivalent model can be obtained as follows
\begin{equation}  \label{33}
\begin{split}
\text{min} & ~~\|Rw-R^I\|^2_2+\lambda\|w\|_{1/2}^{1/2} \\
\text{s.t.}& ~~e^Tw=1 \\
& ~~\eta_i\leq w_{i}\leq \delta_i, \ \ i\in Supp(w) \\
&~~w_{i}=0, \ \ i\notin Supp(w)
\end{split}%
\end{equation}
where $\lambda $ is the regularization parameter, setting $\lambda=\infty$
produces the totally constrained solution (K=0) whereas $\lambda=0$ yields
the unrestricted solution.

Minimization of $L_{1/2}$ constraints is now a widely used technique when
sparse solutions are desirable. In index tracking problem, sparsity also
play a key role in the task of formulating tracking portfolio. In practice,
managers often want to limit the number of assets or the proportion of
investment the tracking problem. Then the index tracking problem can be
regraded as sparse problem. Fortunately, the new model (\ref{33}) can
provide a sparse solution by controlling the parameters $\lambda $.

Furthermore, we define the two indicators to test the consisitency and
out-of-sample prediction error of the sparse index tracking model. The
consisitency of the model is defined by the absolute difference value of the
error between in-sample and the out-of-sample. The smaller consisitency
means the higher consisitency of the index tracking model. Empirical tests
given in section 4 show that our index tracking model has high consisitency,
so it performs well both in-sample and out-of-sample.

Next we give the remark to show the index tracking model based on \bigskip $%
L_{1}$ regularization.

\textbf{Remark 1.} As Brodie et al and his paper say, $L_{1}$ regularization
can be used to solve the index tracking problem with short selling
constraints \cite{Brodie}. without loss of generality, we consider the
no-short selling index tracking problem in this paper. Next we give the
index tracking model by using $L_{1}$ regularization, that is
\begin{equation}
\begin{split}
\text{min}& ~~\Vert Rw-R^{I}\Vert _{2}^{2}+\lambda \Vert w\Vert _{1} \\
\text{s.t.}& ~~e^{T}w=1 \\
& ~~\eta _{i}\leq w_{i}\leq \delta _{i},\ \ i\in Supp(w) \\
& ~~w_{i}=0,\ \ i\notin Supp(w).
\end{split}
\label{l1model}
\end{equation}

\bigskip The next subsection we will give the efficient algortihm to solve
the above model.


\subsection{A hybrid half thresholding algorithm}


As previously discussed, the half thresholding algorithm in the subsection
2.2 is to solve the $L_{1/2}$ regularization without any constrains, but the
model (\ref{33}) is the $L_{1/2}$ regularization with the convex constrains.
In this subsection we propose a hybrid half algorithm to solve the model (%
\ref{33}).

The hybrid half thresholding algorithm is divided two steps, which is to
handle separately the $L_{1/2}$ regularization of selecting the support set
and the quadratic optimization problem that consists in finding the optimal
asset weights for the fixed $K$ stocks. We first consider the unconstrained
case, i.e. the minimization of the objection function of model (\ref{33}),
and then discuss how to deal with the constraints.

In the first step, we discuss the algorithm for minimizing the objective
function of the index tracking model (\ref{33}), that is
\begin{equation}
\text{min}~~\Vert Rw-R^{I}\Vert _{2}^{2}+\lambda \Vert w\Vert _{1/2}^{1/2}.
\label{mubiao}
\end{equation}%
Clearly, the model (\ref{mubiao}) can be regarded as the $L_{1/2}$
regularization if the parameter $A$ and $b$ in $L_{1/2}$ regularization are
replaced by the parameter $R$ and $R^{I}$. Suppose the $w^{n}$ is current
iterate point, an iteration
\begin{equation*}
w_{n+1}=H_{\lambda _{n}\mu _{n},1/2}(w_{n}+\mu _{n}R^{T}(R^{I}-Rw_{n})).
\end{equation*}%
can be naturally defined, which is called half thresholding algorithm for $%
L_{1/2}$ regularization. Furthermore, If we need $K$ stocks to track the
object index, i.e. the model (\ref{mubiao}) can be regarded as the $K$%
-sparsity problem. Incorporated with different parameter-setting strategies
in \cite{54}, the parameters are adopted by
\begin{equation*}
\mu _{n}=\mu _{0},~~\lambda _{n}=\min \{\lambda _{n-1},~~\frac{\sqrt{96}}{9}%
\Vert R\Vert ^{2}\left\vert [B_{\mu _{n}}(w_{n})]_{K+1}\right\vert ^{\frac{3%
}{2}}\},
\end{equation*}%
where
\begin{equation*}
\mu _{0}=\frac{1-\varepsilon }{\Vert R\Vert ^{2}}
\end{equation*}%
with any small $\varepsilon \in (0,1),$ $B_{\mu _{n}}(w_{n})=w_{n}+\mu
_{n}R^{T}(R^{I}-Rw_{n}),$ and $K$. When so doing, an iteration algorithm
will be adaptive, and free from the choice of regularization parameter.

In the second step, through selecting the support set of the tracking
portfolio $w$, we have the optimal asset weights $w_{i}=0$ if $i\notin
Supp(w).$ The nonzero optimal asset weights can be solved by the following
quadratic programming
\begin{equation}
\begin{split}
\text{min}& ~~\Vert \bar{R}w-R^{I}\Vert _{2}^{2} \\
\text{s.t.}& ~~e^{T}w=1 \\
& ~~\eta _{i}\leq w_{i}\leq \delta _{i},\ i\in Supp(w).
\end{split}
\label{3444}
\end{equation}%
where $\bar{R}\in R^{K\times T}$ is the corresponding returns matrix of
stocks with the nonzero weights. There exist very efficient algorithm to
solve the above model(\ref{3444}). In this paper we adopt the Matlab
function (quadprog) to solve it.

\textbf{Remark 2.} Consider the index tracking problem using $L_1$
regularization (\ref{l1model}), we can design the similar algorithm called
to the hybrid half thresholding algorithm. The hybrid LARS algorithm is
divided two parts, First, the Least Angle Regression or LARS \cite{Lars} are
used to solve the following problem,
\begin{equation}  \label{mubiao1}
\text{min} ~~\|Rw-R^I\|^2_2+\lambda\|w\|_{1}.
\end{equation}
The algorithm seeks to solve the above model for a range of value of
regularization parameter $\lambda,$ starting from a very large value, and
gradually decreasing $\lambda$ until the desired value is attainted. As $%
\lambda$ evolves, the optimal solution moves on a piecewise affine path. As
such, to find the needed tracking portfolio with the $K$ nonzero asset
weights. Next the nonzero optimal asset weights can be solved by the same
quadratic programming (\ref{3444}).

\section{\large Empirical results}
In this section, we apply sparse index tracking model (ref:33) and hybrid
half thresholding algorithm described above to conduct optimal tracking
portfolios and evaluate their out-of-sample performance and consistency.

The empirical comparisons are conducted on benchmark problems from the
OR-Library (Beasley. \cite{OR}). For the index tracking problem it contains
the weekly stock prices of the stocks included in major world market
indexes, more specifically, we consider Hang Seng (Hong Kong), DAX 100
(Germany), FTSE (Great Britain), Standard and Poor's 100 (USA), the Nikkei
index (Japan), the Standard and Poor's 500 (USA), Russell 2000 (USA) and
Russell 3000 (USA).

The purpose of empirical tests is to assess the out-of-sample performance of
the sparse index tracking model and hybrid half thresholding algorithm. To
compare performance, two competitive model and algorithm like Torrubiano's
model \cite{Torr} with hybrid optimization approach, $L_1$ model (\ref%
{l1model}) with hybrid LARS approach have been also applied, together with
our model and hybrid half thresholding algorithm. Similar to the experiments
that were carried out by Torrubiano et al \cite{Torr}, the data sets of
weekly returns of the stocks included in the index are partitioned into a
training set containing the first half of the data (145 values) and a test
set with the rest of the data (145 values). The training data sets are used
to find the optimal tracking portfolio, and the testing data sets are used
to estimate the out-of-sample tracking error of the tracking portfolio. The
in-sample and out-of-sample tracking error marked as $TEI$ and $TEO$
respectively. To compare the consistency both in-sample and out-of-sample
and out-of -sample performance of our model with Torrubiano's model \cite
{Torr}and $L_1$ model, we define the following indicators,

$\bullet $\ \ Consistency ($Cons$): This indicator is used to measure the
consistency of the model both in and out of the sample, defined as
\begin{equation*}
Cons=|TEI-TEO|.
\end{equation*}%
Clearly, the smaller value of the Cons means that the model is more stable
both in-sample and out-of-sample.

$\bullet $\ \ Superiority of out-of-sample ($SupO$): We define
\begin{equation*}
SupO=\frac{TEO1-TEO2}{TEO1}\times 100\%,
\end{equation*}%
where $TEO1$ and $TEO2$ are the out-of-sample tracking error of model 1 and
our mode 2. If $SupO>0,$ namely, $TEO2$ is smaller than $TEO1$, i.e. model 2
has the better out-of-sample error than model 1.


The tests were conducted on a personal computer (2.67Ghz, 4Gb of RAM) with
MATLAB 7.9 programming platform (R2009b). The lower and upper bound of the
asset weight set to $\eta_{i}=0.01$ and $\delta_{i}=0.5$, $i=1,2,\cdots,N$.%
\newline

\textbf{A. Comparison with Torrubiano's model }

We present experiments to compare the performance of our model and
Torrubiano's model by using Hang Seng (Hong Kong), DAX 100 (Germany), FTSE
(Great Britain), Standard and Poor's 100 (USA), the Nikkei index (Japan).
The in-sample error and the out-of-sample error of Torrubiano's model are
cited in \cite{Torr} , Results for five data sets are summarized in Table %
\ref{biao41} and Figure 1.

\begin{table}[h!]
{\footnotesize
 \centering\caption{Comparison Torrubiano's model with our model}\label{biao41}
 \renewcommand\arraystretch{1.2}
\begin{center}\begin{tabular}{ccccccccc}
 \hline
Index &Scale & \multicolumn{3}{c}{Our model}   &\multicolumn{3}{c}{Torrubiano's model} & $SupO(\%)$ \\
 &$K$& $TEI2$  &$TEO2$ &$Cons2$  & $TEI1$ & $TEO1$  &$Cons1$ & \\ \hline

Hang &5   &5.81e-5  &4.19e-5  &1.62e-5     &4.14e-5    &7.22e-5   &3.08e-5    &41.91\\
Seng &6   &5.01e-5  &3.85e-5  &1.16e-5     &3.031e-5   &4.76e-5  &1.724e-5   & 19.03\\
$N$=31 &7 &3.56e-5  &2.62e-5  &9.38e-6     &2.37e-5   &3.81e-5   &1.44e-5  &31.15 \\
 &8       &2.61e-5  &2.02e-5  &5.92e-6     &1.91e-5   &2.90e-5   &9.92e-6  &30.36 \\
 &9       &2.31e-5  &1.63e-5  &6.77e-6     &1.62e-5   &2.58e-5    &9.59e-6  &36.85\\
&10       &1.84e-5  &1.64e-5  &2.07e-6     &1.35e-5   &2.06e-5     &7.11e-6 &20.36\vspace{0.2cm} \\

DAX &5     &4.57e-5  &1.20e-4 & 7.40e-5 &2.21e-5 &1.02e-4 &7.97e-5 &\textbf{-17.58}\\
$N$=85 &6  & 3.30e-5 &8.78e-5 &5.47e-5  &1.76e-5 &8.94e-5 &7.17e-5 &1.79 \\
 &7        & 2.41e-5 &9.80e-5 &7.39e-5  & 1.37e-5&8.46e-5 &7.09e-5 &\textbf{-15.83} \\
 &8        & 2.14e-5 &8.97e-5 &6.83e-5 & 1.11e-5& 7.93e-5& 6.82e-5&\textbf{-13.08} \\
 &9        &1.94e-5  &8.80e-5 &6.86e-5 &9.22e-6 & 7.78e-5& 6.85e-5&\textbf{-13.14} \\
&10        &2.96e-5  &2.90e-5&5.68e-5 &8.08e-6 & 7.48e-5&6.67e-5  &61.22  \vspace{0.2cm} \\

FTSE &5    & 1.14e-4 & 9.01e-5& 2.37e-5 & 6.42e-5& 1.58e-4&9.39e-5 & 43.00\\
$N$=89 &6  & 8.30e-5 &8.68e-5 & 3.72e-6 &4.96e-5 &1.12e-4 & 6.23e-5&22.47 \\
     &7    & 7.91e-5 &7.42e-5 &4.87e-6  &3.83e-5 & 9.07e-5& 5.24e-5&18.15 \\
     &8    & 6.24e-5 &6.72e-5 &4.83e-6 &2.90e-5 & 9.66e-5& 6.76e-5&30.45 \\
     &9    & 5.60e-5 & 5.64e-5& 6.19e-6 &2.49e-5 &8.59e-5 & 6.11e-5&34.41\\
     &10   & 4.30e-5 &4.92e-5 &6.19e-6  &2.18e-5 &8.01e-5 &5.82e-5 &38.54  \vspace{0.2cm} \\

S\&P  &5 & 1.21e-4 &1.09e-4 & 1.06e-5    &4.50e-5    &1.14e-4& 6.92e-5& 3.72\\
$N$=98 &6 & 6.80e-5 &8.30e-5 & 1.50e-5   &3.37e-5  &1.01e-4 & 6.70e-5& 17.61\\
     &7 & 8.72e-5 &8.33e-5 & 3.88e-6     & 2.76e-5   &7.80e-5&5.04e-5 &\textbf{-6.80} \\
     &8 &3.89e-5  & 5.98e-5& 2.08e-5     &2.27e-5    &6.76e-5& 4.49e-5&11.66\\
     &9 &7.42e-5  & 4.90e-5& 2.52e-5     & 1.94e-5   &5.91e-5&3.97e-5 &17.05 \\
     &10 &3.99e-5  &4.22e-5& 2.25e-6     &1.66e-5    &5.55e-5&3.89e-5 & 23.96 \vspace{0.2cm} \\

Nikkei &5    &1.26e-4  &1.58e-4 & 3.19e-5 &5.46e-5 & 1.63e-4&1.08e-4 & 2.87\\
$N$=225&6    &1.15e-4  &1.41e-4 & 2.58e-5 & 4.01e-5&1.47-4 &1.07e-4 &3.93 \\
       &7    &8.81e-5 & 1.21e-4& 3.38e-5 &3.36e-5 &1.32e-4 & 9.88e-5&7.93 \\
       &8    &5.94e-5 & 9.34e-5&3.40e-5  &2.60e-5 &1.10e-4 &8.40e-5 &15.08 \\
      &9     &5.96e-5 &8.14e-5 & 2.18e-5  &2.13e-5 & 9.80e-5 &1.68e-5 &17.01\\
      &10    &7.08e-5  &6.96e-5 & 1.29e-6& 1.80e-5 & 6.47e-5& 4.67e-5& \textbf{-7.49}\\ \hline
\end{tabular}\end{center}
}
\end{table}

From the table \ref{biao41} we see that:

$\bullet ~~$ Our model has lower out-of-sample prediction error since $SupO>0
$ at $80\%(=24/30)$. But this is not necessarily the case in the training
sets, the emphasis on our model is to improve the tracking error in testing
data sets for higher prediction effect. We take the FISE index as the
example. In the Figure 1, the in-sample error and the out-of sample error of
our model is in the middle position. i.e. our model has the better out-of
sample error than the Torrubiano's model at the cost of in-sample error.

\begin{figure}[htbp]
  \centering
    \label{fig:subfig:a}
  \includegraphics[width=3.5in]{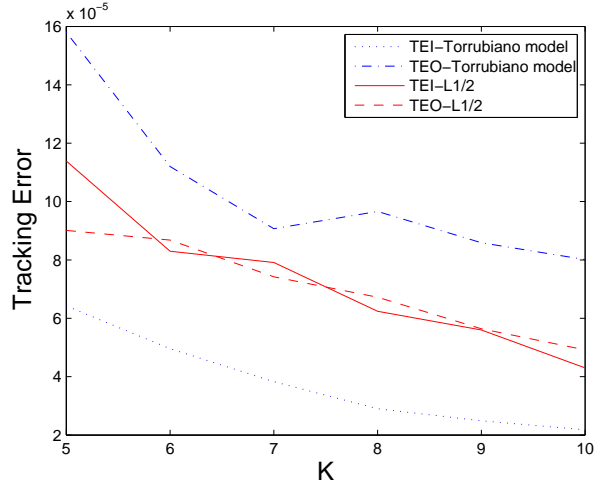}
\caption{Comparsion of Torrubiano model and $L_1/2$ model based on FISE index }
\end{figure}

$\bullet~~$ Our model is more stable than Torrubiano's model since the
consistence indicator $Cons2$ is smaller than $Cons1.$

\textbf{B. Comparison with $L_{1}$ model }

Next we conduct experiments to compare the performance of our model using
hybrid half thresholding algorithm and $L_{1}$ model using hybrid LARS
algorithm. Numerical results are listed in Table $\ref{biao42}$ and Figure 2.

\begin{table}[h!]
{\footnotesize
 \centering\caption{Comparison $L_1$ model with our model}\label{biao42}
 \renewcommand\arraystretch{1.2}
\begin{center}\begin{tabular}{ccccccccc}
 \hline
Index &Scale & \multicolumn{3}{c}{Our model}   &\multicolumn{3}{c}{L1 model} & $SupO(\%)$ \\
 &$K$& $TEI2$  &$TEO2$ &$Cons2$  & $TEI1$ & $TEO1$  &$Cons1$ & \\ \hline

Hang &5   &5.81e-5  &4.19e-5  &1.62e-5     &9.65e-5    &9.35e-5   &2.96e-6    &55.15\\
Seng &6   &5.01e-5  &3.85e-5  &1.16e-5     &5.47e-5   &7.09e-5  &1.62e-5   & 45.73\\
$N$=31 &7 &3.56e-5  &2.62e-5  &9.38e-6     &4.74e-5   &5.64e-5   &9.06e-6  &53.50 \\
 &8       &2.61e-5  &2.02e-5  &5.92e-6     &4.59e-5   &5.95e-5   &1.36e-5  &66.05 \\
 &9       &2.31e-5  &1.63e-5  &6.77e-6     &4.41e-5   &5.13e-5    &7.13e-6  &68.20\\
&10       &1.84e-5  &1.64e-5  &2.07e-6     &4.02e-5   &4.20e-5     &1.79e-6 &61.01\vspace{0.2cm} \\

DAX &5     &4.57e-5  &1.20e-4 & 7.40e-5 &3.26e-5 &1.22e-4 &8.94e-5 &1.86\\
$N$=85 &6  & 3.30e-5 &8.78e-5 &5.47e-5  &2.25e-5 &9.89e-5 &7.64e-5 &11.24\\
 &7        & 2.41e-5 &9.80e-5 &7.39e-5  & 1.66e-5&8.57e-5 &6.90e-5 &\textbf{-14.37} \\
 &8        & 2.14e-5 &8.97e-5 &6.83e-5 & 1.60e-5& 8.04e-5& 6.44e-5&\textbf{-11.48} \\
 &9        &1.94e-5  &8.80e-5 &6.86e-5 &1.49e-5 & 7.81e-5& 6.32e-5&\textbf{-12.62} \\
&10        &2.96e-5  &2.90e-5&5.68e-5 &1.43e-5 & 7.81e-5&6.37e-5  &62.82  \vspace{0.2cm} \\

FTSE &5    & 1.14e-4 & 9.01e-5& 2.37e-5 & 1.06e-5& 1.33e-4& 2.67e-5 &32.29\\
$N$=89 &6  & 8.30e-5 &8.68e-5 & 3.72e-6 &9.94e-5 &1.18e-4 &1.83e-5&26.31\\
     &7    & 7.91e-5 &7.42e-5 &4.87e-6  &8.78e-5 & 1.14e-4& 2.57e-5&34.63 \\
     &8    & 6.24e-5 &6.72e-5 &4.83e-6 &7.61e-5 & 1.16e-4 &4.02e-5&42.19 \\
     &9    & 5.60e-5 & 5.64e-5& 6.19e-6 &5.62e-5 &9.40e-5 &3.41e-5&39.95\\
     &10   & 4.30e-5 &4.92e-5 &6.19e-6  &5.34e-5 &8.75e-5 &3.41e-5 &43.74  \vspace{0.2cm} \\

S\&P  &5 & 1.21e-4 &1.09e-4 & 1.06e-5    &1.01e-4    &1.26e-4& 2.44e-5& 12.39\\
$N$=98 &6 & 6.80e-5 &8.30e-5 & 1.50e-5   &8.15e-5  &9.26e-4 & 1.10e-5& 10.36\\
     &7 & 8.72e-5 &8.33e-5 & 3.88e-6     &5.56e-5   &7.51e-5&1.95e-5 &\textbf{-10.95} \\
     &8 &3.89e-5  & 5.98e-5& 2.08e-5     &4.44e-5    &6.80e-5& 2.36e-5&12.13\\
     &9 &7.42e-5  & 4.90e-5& 2.52e-5     &4.27e-5   &5.98e-5&1.71e-5 &18.00 \\
     &10 &3.99e-5  &4.22e-5& 2.25e-6     &4.22e-5    &5.73e-5&1.51e-5 & 26.39 \vspace{0.2cm} \\

Nikkei &5    &1.26e-4  &1.58e-4 & 3.19e-5 &1.48e-5 & 2.10e-4&6.24e-5 & 24.72\\
$N$=225&6    &1.15e-4  &1.41e-4 & 2.58e-5 & 1.31e-5&2.20-4 &8.93e-5 &35.87 \\
       &7    &8.81e-5 & 1.21e-4& 3.38e-5 &1.18e-5 &1.82e-4 & 6.39e-5&32.85 \\
       &8    &5.94e-5 & 9.34e-5&3.40e-5  &1.08e-5 &1.66e-4 &5.83e-5 &43.69 \\
      &9     &5.96e-5 &8.14e-5 & 2.18e-5  &9.89e-5 &1.62e-4 &6.29e-5 &49.92\\
      &10    &7.08e-5  &6.96e-5 & 1.29e-6& 9.47e-5 & 1.59e-4& 6.42e-5& 56.25\\ \hline
\end{tabular}\end{center}
}
\end{table}

From the Table 2 we see that our model has lower out-of-sample prediction
error than $L_{1}$ model since $SupO>0$ at $87\%(=26/30)$. Moreover, we find
our model can provide more sparse solution to track the object index. The
Figure 2 shows this results. The out-of-sample prediction error of $L_{1}$
model using $K=10$ stocks is the same to our model using $K=5$ stocks.

\begin{figure}[htbp]
  \centering
    \label{fig:subfig:a} 
  \includegraphics[width=3.5in]{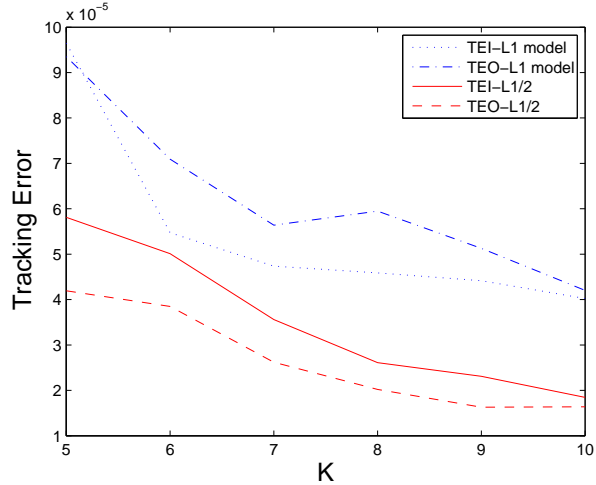}
\caption{Comparsion of $L_1$ model and $L_1/2$ model based on Hang Seng index }
\end{figure}

Finally, we discuss the large scale index tracking problem, i.e., Standard
and Poor's 500($N=457$), Russell 2000 ($N=1318$) and Russell 3000($N=2151$).
Since the number of stocks included in indexes is very large, we select the
number of the tracking stocks $K=10,20,30,40,50$, the numerical results are
listed in Table \ref{biao43}. According to the Table \ref{biao43}, the value
$SupO>0$ for all cases, it is shown that our model has better out-of-sample
prediction error than $L_1$ model.

\begin{table}[h!]
{\footnotesize
 \centering\caption{Comparison $L_1$ model with our model}\label{biao43}
 \renewcommand\arraystretch{1.2}
\begin{center}\begin{tabular}{ccccccccc}
 \hline
Index &Scale & \multicolumn{3}{c}{Our model}   &\multicolumn{3}{c}{L1 model} & $SupO(\%)$ \\
 &$K$& $TEI2$  &$TEO2$ &$Cons2$  & $TEI1$ & $TEO1$  &$Cons1$ & \\ \hline

S\&P   &10   &1.49e-4  &4.70e-4  &2.92e-4     &1.08e-4    &3.44e-4   &2.36e-4    &26.88\\
 &20         &7.93e-5  &2.76e-4  &1.97e-4     &3.27e-5   &1.73e-4  &1.41e-4   & 37.18\\
$N$=457 &30  &5.03e-5  &2.23e-4  &1.72e-4     &3.78e-5   &1.61e-4   &1.23e-4  &27.69\\
 &40       &3.42e-5  &1.68e-4  &1.34e-4     &3.81e-5   &1.10e-4   &7.17e-5  &34.78 \\
 &50        &2.57e-5  &1.42e-4  &1.16e-4     &4.18e-5   &1.14e-4    &7.27e-5  &19.17\vspace{0.2cm} \\

Russell &10     &4.62e-4  &5.98e-4 & 1.37e-4 &2.29e-4 &5.77e-4 &3.48e-4 &3.52\\
$N$=1318 &20  & 1.56e-4 &4.34e-4 &2.78e-4  &1.20e-4 &3.83e-4 &2.63e-4 &11.86\\
 &30        &1.06e-4 &4.08e-4 &3.03e-4  & 1.30e-4&3.23e-4 &1.93e-4 &20.94 \\
 &40        & 5.48e-5 &3.20e-4 &2.66e-4 & 7.89e-5& 2.32e-4& 1.53e-4&27.68 \\
 &50        &5.22e-5  &2.89e-4 &2.36e-4 &9.78e-5 & 2.62e-4& 1.65e-4&9.06  \vspace{0.2cm} \\

Russell &10    & 1.26e-4 & 4.91e-4& 3.65e-4 &3.78e-4& 3.97e-4& 1.94e-5 &19.14\\
$N$=2151 &20  & 7.44e-5 &3.09e-4 & 2.34e-4  &1.22e-4 &2.37e-4 &1.15e-4&23.26\\
     &30    & 3.88e-5 &2.37e-4 &1.98e-4     &1.27e-4 & 2.28e-4& 1.00e-4&3.96 \\
     &40    & 3.39e-5 &2.07e-4 &1.73e-4     &8.45e-5 & 1.06e-4 &1.22e-4&0.44 \\
     &50    & 3.71e-5 & 1.69e-4& 1.32e-4    &1.31e-4 &1.67e-4 &3.56e-5&1.44\vspace{0.2cm} \\
\hline
\end{tabular}\end{center}
}
\end{table}

\section{\large Conclusions}

Index tracking is a passive financial strategy that
aims at replicating the performance and risk-profile of a given index.
One of the most common approaches to tackle the index tracking
problem consists of minimizing a given tracking error measure while limiting
the maximum number of assets held in the portfolio. Having few active
positions reduces the administrative and transaction costs and avoids detaining
very small and illiquid positions, especially when the index has a large
number of constituents. However, imposing an upper bound on the number
of constituents of the tracking portfolio makes the optimization problem NP-Hard.
 Different quantitative approaches have
been proposed to tackle such an optimization problem. Most approaches rely on
search heuristics. On the other hand, $L_1$ regularization methods have found application in mean-variance portfolio settings in order to promote the
identification of sparse portfolios with
good out-of-sample properties and low turnover. However, the  $L_1$ regularization
approach is ineffective in index tracking problem, since the index tracking problem has budget and
no-short selling constraints.

In this paper We have used a new constrains $\|w\|_{1/2}^{1/2}=K$ of tracking
portfolio's weight to replace the cardinality constrains
$\|w\|_{0}=K$ which equals to $\sum_{i=1}^N Z_i=K,~Z_i=0~or~1$. A
new sparse index tracking model was established by minimizing
tracking error. Different to the other models of stock index
tracking, our model has high consistency and out-of-sample
prediction error, it also can preserve sparsity of the optimal tracking
portfolio as much as possible. Meanwhile, since the half threshoding algorithm is the fast solver of $L_{1/2}$ regularization, we have
extended the half threshoding algorithm to hybrid half thresholding algorithm for solving
the proposed index tracking model. The algorithm is fast and efficient with appropriate parameters selection
for the sparse index tracking model. Furthermore, it is simple, very convenient in use, and can be applied to large scale problems.

We have tested performance of our model and algorithm on the eight
data sets from OR-library. Numerical results have shown that the
sparse index tracking model and hybrid half thresholding
algorithm have high consisitency and better out-of-sample prediction
ability. We believe the sparse index tracking model and projected
half thresholding algorithm can provide useful reference to the
manager of index derivatives. Next we plan to extend our results to
the index tracking problems with transaction costs.

\end{document}